\documentclass{article}
\usepackage{amssymb}
\usepackage{amsmath}
\usepackage[T1]{fontenc}

\begin{document}
\newcommand{\ov}{\overline}
 \newcommand{\un}{\underline}
\newcommand{\proof}{\bf {Proof:} \rm}
\newcommand{\dirac}{ {\bf D}}
\newcommand{\cauchyr}{ {\mathcal D}_{h}}
\newcommand{\diracc}{ \tilde{D}_{h}}
\newcommand{\difsc}{ d_{h}}
\newcommand{\dif}{ \partial_{h}}
\newcommand{\lapl}{ \Delta_{h}}
\newcommand{\laplc}{ \tilde{\Delta}_{h}}
\newcommand{\euler}{ E_{h}^{\pm}}
\newcommand{\gam}{ \Gamma_{h}^{\pm}}
\newcommand{\K}{ {\bf K}}

\newcommand{\Grad}{{\rm grad}}
\newcommand{\Curl}{{\rm curl}}
\newcommand{\Div}{{\rm div}}
\newcommand{\Sc}{{\rm Sc\,}}
\newcommand{\wzwo} {\stackrel{\circ}{{\cal W}^1_2}}
\newcommand{\wzwoz} {\stackrel{\circ}{{\cal W}^2_2}}
\newcommand{\wzwok}{\stackrel{\circ}{\cW_2^k}}
\newcommand{\wzwokl}{\stackrel{\circ}{\cW_2^{k,l}}}

\newcommand{\BR}{{\mathbb R}}
\newcommand{\BC}{{\mathbb C}}
\newcommand{\BN}{{\mathbb N}}
\newcommand{\BZ}{{\mathbb Z}}

\newcommand{\x}{\underline{x}}
\newcommand{\m}{\underline{m}}
\newcommand{\p}{\underline{p}}
\newcommand{\q}{\mathfrak{q}}
\newcommand{\f}{\mathfrak{f}}
\newcommand{\edge}{\mathfrak{e}}
\newcommand{\g}{\mathfrak{g}}
\newcommand{\G}{{\bf G}}

\newcommand{\dext}{{\bf d}}
\newcommand{\e}{{\bf e}}
\newcommand{\vv}{{\bf v}}
\newcommand{\bb}{{\bf b}}
\newcommand{\C}{{\bf c}}

\newcommand{\ob}{{\bf b}}
\newcommand{\cW}{{W}}
\newcommand{\cL}{{L}}
\newcommand{\cC}{{C}}
\newcommand{\cM}{{M}}

\newcommand{\cl}{C \kern -0.1em \ell}

\newcommand{\qed}{$\blacksquare$}
\newtheorem{theorem}{Theorem}[section]
\newtheorem{remark}{Remark}[section]
\newtheorem{lemma}{Lemma}[section]
\newtheorem{proposition}{Proposition}[section]
\newtheorem{corollary}{Corollary}[section]
\newtheorem{definition}{Definition}[section]
\newtheorem{example}{Example}[section]

\title{On a correspondence principle between discrete differential forms, graph structure and multi-vector calculus on symmetric lattices \footnote{Dedicated to the memory of Jarol\'{i}m Bur\v{e}s} }

\author{Nelson Faustino$^\dagger$,~~Uwe K\"ahler\thanks{Departamento de Matem\'atica, Universidade de
Aveiro, Portugal} \\ \ \\ {\small \{nfaust,ukaehler\}@ua.pt }} \maketitle

\begin{quote}
\begin{abstract}
Based on \cite{DH94}, we introduce a bijective correspondence between first order
differential calculi and the graph structure
of the symmetric lattice that allows one to encode completely the interconnection structure of the graph
in the exterior derivative. As a result, we obtain the Grassmannian character of the lattice as well as the mutual
commutativity between basic vector-fields on the tangent space.

This in turn gives several similarities between the Clifford setting and the algebra of endomorphisms
endowed by the graph structure, such as the hermitian structure of the lattice as well as the
Clifford-like algebra of operators acting on the lattice. This naturally leads to
a discrete version of Clifford Analysis.
\end{abstract}
\end{quote}


{\bf keywords}: Discrete differential forms, symmetric lattice reduction, hermitian structure.

\section{Introduction}
There are many reasons for studying discrete structures in
mathematics and physics.

As it was shown in a series of papers~\cite{GH03,FGHK06,CFV07}, one
of these reasons is the numerical treatment of problems related to
potential theory and boundary values problems, where the development
of discrete theoretical counterparts of a continuous operator
calculus leads to well-adapted numerical methods. Another reason is that lattices in physics exhibit some
noncommutative geometric nature and correspond to a regularization
process in field theory as well as a naive approximation of the
topology of space or space-time, that may at high energies exhibit a
topology different from that of the continuum.

Although noncommutative geometry has deep roots into quantum mechanics, this notion
has been introduced by A.~Connes in his extension of the calculus of differential forms and the
de Rham homology of currents \cite{Connes86}.

Among the first implications was the construction of a
Yang-Mills-Higgs theory employing the $\mathbb{C}(X)\oplus\BC(X)$
and later $(\mathbb{H} \oplus \BC) \otimes \BC(X)$ algebra, where
$\BC$ and $\mathbb{H}$ denote the field of complex numbers and the skew-field of quaternions,
respectively, and $\BC(X)$ is the algebra of functions over
$X$~\cite{CL90}. Since then, there exists a growing interest among
theorists in studying noncommutative geometry. A major reason for this is that noncommutative geometry is
ultimately related to quantum groups~\cite{Schupp93}. The latter are
connected with some important aspects of physics, such as quantum
spin chains, conformal field theories, quantum integrable models and
so on.

There also exists another proposal by A.~Dimakis et al.~\cite{DHS93}
in which the coordinates are kept commutative but there exists
non-commutativity between the coordinates and their differentials.
In this particular case the continuity is lost and the space
acquires a canonical lattice structure.

Regardless of the last viewpoint one looks to discrete structures as
a certain kind of differential calculi over discrete sets. First
order differential calculi on discrete sets were found to be in
bijective correspondence with graph structures~\cite{DH94}, where
the vertices of a graph are given by elements of the set and neither
multiple edges nor loops are admitted. In particular, this supplies
the elements of the set with neighborhood relations.

The power of differential calculus rests mostly on its intrinsic
character and on the algebraic structure which is the Grassmann
algebra. In spite of their powerful nature, Grassmann algebras do not
incorporate some concepts that are crucial in physics like the
concept of a spinor and of a Dirac operator. This is one among many
reasons why the correspondence between Clifford algebras and lattice
structure should be established {\it a priori}.

It is well known that Clifford algebras can be defined to lattice
structures in terms of the cup algebra of simplicial homology
theory. The idea was therefore to start from the algebra of
endomorphisms of the vector space of cochains, and look for a
natural algebraic structure in this space just like for Clifford
algebras in the continuous case~\cite{Vaz97}. This approach, however, has some limitations, particularly in
relation to lattice gauge theories, where the forward/backward
differences $\partial_h^{\pm j}$ are replaced by the symmetric
differences $\frac{1}{2}(\partial_h^{- j}+\partial_h^{+ j})$. In
particular, we want to write a lattice version of the Dirac operator
which splits a lattice version of the Laplace operator.

One way to overcome this problem was proposed by Wilson
in~\cite{Wilson74}, by adding an extra term to the lattice
version of the Dirac operator such that an extra fermion acquires
mass and of a order of the cut-off. Another way was recently
proposed by Faustino, K\"ahler, and Sommen in~\cite{FKS07} using the
splitting of the standard Clifford basis vectors $\e_j$ in
$\e_j=\e_j^-+\e_j^+$ such that Dirac operators on lattices are
constructed by using superpositions of the type
$\e_j^+\partial_h^{+j}+\e_j^-\partial_h^{-j}$.

Although both approaches seem to be very promising tools, however,
they have some limitations, particularly in relation to the lattice
structure. The first in a certain sense requires the usage of
second-order operators to define a Dirac operator while the second
at first seems to be quite artificial. We will show that the second
approach comes natural as a special case from the hermitian
structure of the lattice. In our opinion, a better geometric
understanding about the concept of Dirac operators is crucial to the
formulation of well-adapted theoretical approaches on lattices,
namely the theory of discrete monogenic functions.

This paper is organized as follows: first we will introduce the
discrete differential geometry setting needed to describe the
notions of exterior and interior product in terms of differential
forms and dual connections between differential forms and vector
fields on the lattice, respectively. Such an axiomatization can be
found in~\cite{DH03} and was inspired by the works of Vaz
\cite{Vaz97},~Sommen~\cite{Sommen97} and Brackx,~Delanghe,~Sommen~\cite{BDS05}.

Next we associate an un-oriented graph with a differential calculus
on $\BZ^n$ by truncating most of the non-local links except the
nearest neighboring links (i.e. a symmetric lattice reduction). This
will be fundamental to prove some well-known properties of Grassmann
algebras and vector-fields, namely the exterior product rule of
differential forms and the mutual commutativity between vector
fields.

Afterwards we introduce the necessary endomorphisms such as
the interior and exterior product which enables us to establish some
similarities between the symmetric structure of the lattice and
complex Clifford algebras $\BC_{2n}$. The resulting Clifford structure
will be compared with the Clifford-like structure on the symmetric lattice
introduced by Kanamori and Kawamoto in~\cite{KK04}.

On the last section, we introduce the Dirac operators and the
corresponding vector variables on symmetric lattices using the
hermitian setting (see, for instance,~\cite{SS02},~\cite{Bures07}
and~\cite{BSS07} and references given therein). We do not claim that
the obtained construction of discrete Dirac operators is new. It can
also be found in the works of Vaz~\cite{Vaz97}, Kanamori,
Kawamoto~\cite{KK04} and Forgy, Schreiber~\cite{FS04}.

We also establish some basic intertwining relations between these operators
which allow us to derive the Euler operator and also the discrete
counterpart of homogeneous polynomials on the symmetric lattice.
The resulting notion of discrete homogeneity will be compared with
the hermitian homogeneity introduced in~\cite{BSS07}.

\section{Discrete differential geometry}

In this section, we present a short resume of the well-known notions
of discrete differential geometry following the works of Dimakis and
M\"uller-Hoissen presented in \cite{DH94,DH97,DH03}. The same
approach can also be found in~\cite{KK04} as well as in~\cite{FS04}.

\subsection{Universal Differential Algebra on a Lattice}
We consider lattice functions as maps from a lattice point $l$ to a
complex number. These functions jointly with pointwise addition and
multiplication constitute an unital algebra $\mathcal{A}$ which is
associative and abelian, and possesses unity ${\bf 1} \in
\mathcal{A}$.

The universal differential calculus corresponds to the pair $(\dext,\Lambda^*\mathcal{A})$, where $\Lambda^*\mathcal{A}$ is a $\BZ-$graded associative algebra (over $\BC$)
\begin{eqnarray*}
\Lambda^*\mathcal{A}=\sum_{r \geq 0}\Lambda^r\mathcal{A}
\end{eqnarray*}
where $\Lambda^r\mathcal{A}$ are $\mathcal{A}-$bimodules, (i.e. can be multiplied from the left and right by elements of $\mathcal{A}$), and $\dext:\Lambda^r\mathcal{A} \rightarrow \Lambda^{r+1}\mathcal{A}$ a $\BC-$linear map which is nilpotent and obeys the graded Leibniz rule:
\begin{eqnarray}
\label{nilp}\dext(\dext \omega_r)=0, \\
\label{leibniz}\dext(\omega_r \omega_s)=(\dext \omega_r)\omega_s+(-1)^r\omega_r
\dext \omega_s,
\end{eqnarray}
where $\omega_r\in \Lambda^r\mathcal{A}$ and $\omega\in \Lambda^*\mathcal{A}$.

From these definitions it can be easily seen that
\begin{eqnarray}
\label{dext1}\dext {\bf 1}={\bf 0}
\end{eqnarray}

Assuming that $\Lambda^0\mathcal{A}:=\mathcal{A}$, we will
furthermore exclusively consider differential algebras $\Lambda^*
\mathcal{A}$, i.e. we require that $\dext$ generates the spaces
$\Lambda^r\mathcal{A}$.

Let $\mathcal{L}$ be a denumerable set and $f$ a $\BC-$valued
functions over $\mathcal{L}$:
\begin{eqnarray*}
f:\mathcal{L} \rightarrow \BC, & l \mapsto f(l)=f_l.
\end{eqnarray*}

Then the algebra is generated by the set of discrete delta-functions
$\{ \bb_{l}\}_{l}$
$$ \bb_l(m)=\delta_{lm},$$
where $\delta_{lm}$ is the standard Kronecker symbol.

Indeed, any function $f \in \mathcal{A}$ can be expanded as $f=\sum_l f_l \bb_l.$
This also gives rise to the following properties:
\begin{eqnarray*}
\bb_l \bb_m=\delta_{lm}\bb_l, \\
\sum_{l}\bb_l={\bf 1},
\end{eqnarray*}
which reflects the pointwise product of functions and assures the completeness of the basis.
From the above properties we obtain
\begin{eqnarray}
\label{prop2} \bb_l \dext \bb_{m}=-\dext \bb_{l}\bb_m+\delta_{lm}\dext \bb_m \\
\label{prop3}\sum_{l}\dext \bb_l={\bf 0}
\end{eqnarray}
by using (\ref{leibniz}) and (\ref{dext1}).

Using these discrete delta functions $\Lambda^r\mathcal{A}$ can be
written in terms of the elements
\begin{eqnarray}
\label{prop4}\bb_{m^0,m^1,\ldots,m^{r-1},m^r}:=\left\{
\begin{array}{ccc}
\bb_{m^0} \dext \bb_{m^1} \ldots \dext \bb_{m^{r-1}}\dext \bb_{m^{r}} ,& m^0 \neq m^1,\ldots, ~m^{r-1} \neq m^{r}\\
0,& \mbox{otherwise}
\end{array}
\right.
\end{eqnarray}
From (\ref{prop2}) we obtain,
\begin{eqnarray*}
\bb_{m^0,m^1,\ldots,m^{r-1},m^r}=\bb_{m^0}\bb_{m^1,\ldots,m^{r-1}}\bb_{m^r}
=\bb_{m^0}\dext \bb_{m^1} \ldots \dext \bb_{m^{r}}\bb_{m^{r}}
\end{eqnarray*}

Therefore, all the elements $\bb_{m^0,m^1,\ldots,m^{r-1},m^r}$ are built by concatenating elements of the form $\bb_{m^j,m^{j+1}}$:
\begin{eqnarray}
\label{simplex}\bb_{m^0,m^1,\ldots,m^{r-1},m^r}=\bb_{m^0,m^1}\ldots\bb_{m^{r-1},m^r}.
\end{eqnarray}

This naturally leads to the $\dext-$action
\begin{eqnarray}
\label{prop5}\dext \bb_{m^0,m^1,\ldots,m^{r-1},m^r}= \sum_{l\in
\mathcal{L}} \sum_{s=0}^r (-1)^s \bb_{m^0,m^1,\ldots,m^{s-1},l,m^s,
\ldots ,m^r},
\end{eqnarray}
which represents the action of $\dext$ on an $r-$form.

All this shows that $\bb_{m^0,m^1,\ldots,m^{r-1},m^r}$ describe
simplicial $r-$paths which are related to the connectivity of the
discrete space and assures the nil-potency of the discrete
differential operator $\dext$, when acting on (\ref{simplex})
\begin{eqnarray*}
\dext(\dext \bb_{m^0,m^1,\ldots,m^{r-1},m^r})=0
\end{eqnarray*}

In particular, every 1-path, described in terms of non-vanishing
elements $\bb_{l,m}$, corresponds to an edge connecting nodes $l$
and $m$. This induces a graph structure on the set $\mathcal{L}$.

In order to handle the interconnection of points of $\mathcal{L}$,
it is helpful to consider some additional structure: To make use of
the additive structure of the lattice, it is sufficient to assume
that the set $\mathcal{L}$ is equipped with a group addition
$l+m=p$, which allows to {\it go from $l$ in the direction $m$ and
arrive at p} on the discrete space and a group inverse $-m$ defined
by $p+(-m)=l$ which corresponds to the direction $m$ oriented in the
opposite side. This particular class of graphs are Cayley graphs,
i.e. discrete groups.

The left action of a Cayley graph on the lattice will be denoted by
the translation operator $T_l$ defined by
\begin{eqnarray}
T_l f=\sum_m f_{m+l}\bb_m&=&\sum_m f_{m}\bb_{l-m}.
\end{eqnarray}
Then every node $l \in \mathcal{L}$ gives rise to the $1-$forms
\begin{eqnarray}
\label{1formdiff}\Theta^l= \sum_m \bb_{m,m+l}, & \mbox{if}~\sum_{m}
\Theta^m \bb_l \neq 0,
\end{eqnarray}
which are obviously left invariant in the following sense
\begin{eqnarray*}
\label{shifti} T_p \Theta^l=0 &\Longleftrightarrow& \sum_m \bb_{m,m+l}=\sum_m \bb_{m+p,m+p+l}.
\end{eqnarray*}
A further interesting property of the forms $\Theta^l$ (and indeed
of the discrete calculus as a whole) is that they do not commute
with functions, but instead induce translations on them:
\begin{eqnarray}
\label{transl}
\begin{array}{ccc} \bb_l \Theta^m, &=&\Theta^m \bb_{l+m}
\\ \Theta^l f&=&(T_l f)\Theta^l.
\end{array}
\end{eqnarray}

Using this quantity we can rewrite (\ref{simplex}) as
\begin{eqnarray*}
\bb_{m^0,m^1,\ldots,m^{r-1},m^r}=\bb_{m^0}\Theta^{m^1-m^0}\ldots\Theta^{m^r-m^{r-1}}
\end{eqnarray*}
In terms of 1-forms $\Theta^{m^j}$ every element $\omega_r\in
\Lambda^r\mathcal{A}$ is given by
\begin{eqnarray}
\label{rform}\omega_r= \sum_{j=0}^r \sum_{m^j} F_{m^0,m^1,\ldots,m^r}\Theta^{m^1-m^0}\ldots\Theta^{m^r-m^{r-1}}
\end{eqnarray}
where
\begin{eqnarray}\label{function}F_{m^0,m^1,\ldots,m^r}= f_{m^0,m^1,\ldots,m^r}
\bb_{m^0}\in \mathcal{A}.\end{eqnarray}

Let us finally remark that the 1-forms (\ref{1formdiff}) allow to
write the exterior derivative in terms of the the suggestive action
\begin{eqnarray}
\label{diffun}\dext f =\sum_{l,m}f_l
(\bb_{m,l}-\bb_{l,m})=\sum_{l}(T_l f-f)\Theta^l.
\end{eqnarray}

\subsection{Vector fields and dual connections on a lattice}
Let $\mathcal{T}$ denote the dual space of $\Lambda^1\mathcal{A}$ as
a complex vector space. Let $\{\partial_{m,l}\}$ be the basis
in $\mathcal{T}$, dual to the basis $\{ \bb_{l,m}\}$, i.e.
\begin{eqnarray}
\label{contraction}\langle \bb_{l,m} ,\partial_{p,q}
\rangle_0=\delta_{l,q}\delta_{m,p}.
\end{eqnarray}
$\mathcal{T}$ is turned into a corresponding (left or right)
$\mathcal{A}$-module by introducing the left and right actions
\begin{eqnarray}
\label{contraction2}\langle \alpha , f \cdot X \rangle_0=\langle \alpha f ,X \rangle_0,
& \langle \alpha , X \cdot f \rangle_0=\langle f \alpha ,X
\rangle_0.
\end{eqnarray}
As a consequence we have
\begin{eqnarray*}\bb_l \cdot
\partial_{p,m}=\delta_{m,p}
\partial_{p,m}, & \partial_{p,m} \cdot
\bb_l=\delta_{p,l}\partial_{p,m}.
\end{eqnarray*}
Thus, any element $X \in \mathcal{T}$ can be uniquely decomposed
into
$$X=\sum_{l,m}X_l^m \partial_{m,l},$$
where the summation $\sum_{l,m}$ runs over all $l,m \in \mathcal{L}$
for which there is an edge from $l$ to $m$ associated to $\Lambda^1
\mathcal{A}$.

Now we introduce a duality contraction $\langle \cdot , \cdot
\rangle$ on $\Lambda^1 \mathcal{A}$ as a right $\mathcal{A}-$module
and $\mathcal{T}$ as a left $\mathcal{A}-$module by putting
\begin{eqnarray}
\langle \bb_{l,m},X\rangle:=\bb_{l}\langle \bb_{l,m},X\rangle_0.
\end{eqnarray}
Then we have
\begin{eqnarray}
\label{contraction3}
\langle f \alpha , X \cdot h \rangle = f \langle \alpha,X\rangle h,
& \langle \alpha, f \cdot X \rangle= \langle \alpha f ,X \rangle.
\end{eqnarray}
Moreover, the elements of $\mathcal{T}$ become operators on
$\mathcal{A}$ {\it viz}
$$ X(f):=\langle \dext f,X \rangle . $$
Using the Leibniz rule for $\dext$, one proves
\begin{eqnarray*}
X(fh)=f X(h)+(h \cdot X)(f).
\end{eqnarray*}
The space
\begin{eqnarray*}
\mathcal{T}_l:=\{ X \cdot \bb_l~:~X \in \mathcal{T} \}
\end{eqnarray*}
may be regarded as the tangent space at $l \in \mathcal{M}$, dual to
the subspace
\begin{eqnarray*}
\label{cotangent}\Lambda_l^1 \mathcal{A}:=\bb_l \Lambda^1
\mathcal{A}=\{ \bb_l f~:~f \in \Lambda^1\mathcal{A}\}
\end{eqnarray*}
 with respect to the duality contraction
$\langle \cdot, \cdot \rangle$. The space $\Lambda_l^1 \mathcal{A}$ may be regarded as the cotangent
space at $l \in \mathcal{M}$.

As a consequence of this construction, $\{\partial_{m,l} \}_m$ is a
basis of $\mathcal{T}_l$ which is dual to the basis
$\{\bb_{l,m}\}_m$ of $\Lambda^1_l \mathcal{A}$.

 We thus have the following decompositions
\begin{eqnarray*}
\Lambda^1 \mathcal{A}&=&\sum_{l \in
\mathcal{M}}\bigoplus\Lambda^1_l \mathcal{A}, \\ \mathcal{T}&=&\sum_{l \in
\mathcal{M}}\bigoplus \mathcal{T}_l.
\end{eqnarray*}

\section{Hermitian and lattice structure}

\subsection{Discrete differential calculi and graph structure} \label{graphS}
In the last section we introduce universal differential algebra
which is the lattice counterpart of a differential structure on the
ordinary cases. However is sometimes relatively far away from the
concrete lattice formulations, since the the construction of
$\bb_{m,l}$ is lacking, at least initially, {\it the neighboring
structure}.

In order to obtain an appropriate lattice formulation, we need to
make a reduction of edges by truncating most of the non-local except
for the nearest neighboring nodes. Indeed this can be obtained from the first order differential calculus as a quotient $\Lambda^1\mathcal{A}/\mathcal{J}^1$, where $\mathcal{J}^1$ is a submodule of $\Lambda^1\mathcal{A}$ generated by non-vanishing elements of the form $\bb_{\underline{m},\underline{p}}$.

This kind of differential calculi appears frequently on literature. In particular, it seems that most bicovariant differential calculi on quantum groups have this property, see e.g \cite{HR93}. There are also several examples reduced first order differential calculi on discrete groups and cellular networks (see e.g. \cite{DH94,DH03,Req98} and the references given there).

In our case, we will consider the type of reduction proposed by Dimakis and
M\"uller-Hoissen in~\cite{DH94} for the symmetric lattice.

From now on, we will consider the case of our surface $\mathcal{M}$
being the Euclidean space $\BR^n$ partitioned in $n-$dimensional
simplicial complexes, whose vertices form a $n-$dimensional lattice
$\mathcal{L}$, isomorphic to $\BZ^n$.

We start to consider a simplicial complex labeled by the set of
vectors $ \{ \vv_{j}:~j\in [n] \cup [n]' \}$, where $[n]=\{
1,2,\ldots,n\},~[n]'=\{ 1',2',\ldots,n'\}$, $':[n] \cup [n]'
\rightarrow [n] \cup [n]'$ is an involutive permutation which maps
$j \mapsto j'$ and $j' \mapsto j$ and $j \mapsto \vv_{j}$, $j'
\mapsto -\vv_{j}$

This lead to the lattice reduction
\begin{eqnarray}
\label{latticeR}\bb_{\m,\p}~\left\{
\begin{array}{cc}
\neq 0 & \mbox{if}~\p=\m+\vv_{j},~j \in [n] \cup [n]' \\
=0 & \mbox{otherwise}
\end{array}
\right. & (\m,\p \in \mathcal{L}).
\end{eqnarray}
This reduction implies that there are $2n$
nonzero differentials $\Theta^{\p}$ with
\begin{eqnarray}
\Theta^{\vv_{j}}=\sum_{\m \in \mathcal{L}}\bb_{\m,\m+\vv_{j}}, & j
\in [n] \cup [n]'.
\end{eqnarray}
Let us remark that, by construction, the vectors $\vv_{j}$ and
$\vv_{j'}$ have equal lengths. Then the translations $\m+\vv_{j}$
and $\m+\vv_{j'}$ are {\it symmetric} to each other with respect to
the hyperplane $\mathcal{H}$ which contains the point $\m$ of the
coordinate system in $\BR^n$ and which is perpendicular to the edge
which links the nodes $\m+\vv_{j}$ and $\m+\vv_{j'}$.

This means that in our lattice reduction, the edges between the nodes
$\m,\m+\vv_{j}$ and $\m,\m+\vv_{j'}$ are being kept and thus links
are {\it symmetric} with respect to their orientation.

Therefore, the important automorphisms on $\Lambda^*\mathcal{A}$ leaving the
symmetric lattice reduction invariant are the main involution
defined by
\begin{eqnarray}
\label{involution}\begin{array}{ccc}
(\omega_r\omega_r)'=\omega_r'\omega_s', \\
(F_{\vv_{j_1},\ldots,\vv_{j_r}}\Theta^{\vv_{j_1}}\ldots\Theta^{\vv_{j_r}})'=F_{j_1,\ldots,j_r}(\Theta^{\vv_{j_1}})'\ldots(\Theta^{\vv_{j_r}})',
\\
(\Theta^{\vv_{j}})'=\Theta^{\vv_{j'}},
\end{array}
\end{eqnarray}
the reversion given by
\begin{eqnarray}
\label{reversion}\begin{array}{ccc}
(\omega_r\omega_r)^{\widetilde{~}}=\omega_s^{\widetilde{~}}\omega_r^{\widetilde{~}}, \\
(F_{\vv_{j_1},\ldots,\vv_{j_r}}\Theta^{\vv_{j_1}}\ldots\Theta^{\vv_{j_r}})^{\widetilde{~}}=(\Theta^{\vv_{j_r}})^{\widetilde{~}}\ldots(\Theta^{\vv_{j_1}})^{\widetilde{~}}F_{\vv_{j_1},\ldots,\vv_{j_r}},
\\
(\Theta^{\vv_{j}})^{\widetilde{~}}=-\Theta^{\vv_{j'}},
\end{array}
\end{eqnarray}
and, finally, the $\dag-$~conjugation defined by
\begin{eqnarray}
\label{hconjugation}\begin{array}{c}(\omega_r\omega_r)^{\dag}=\omega_s^{\dag}~\omega_r^{\dag}, \\
(F_{\vv_{j_1},\ldots,\vv_{j_r}}\Theta^{\vv_{j_1}}\ldots\Theta^{\vv_{j_r}})^{\dag}=(\Theta^{\vv_{j_r}})^{\dag}\ldots(\Theta^{\vv_{j_1}})^{\dag}\overline{F}_{\vv_{j_1},\ldots,\vv_{j_r}},
\\
(\Theta^{\vv_{j}})^{\dag}=-\Theta^{\vv_{j'}}, \end{array}
\end{eqnarray}
(we remark that
$\overline{F}_{\vv_{j_1},\ldots,\vv_{j_r}}=\overline{f}_{\vv_{j_1},\ldots,\vv_{j_r}}
\bb_{\m}$ stands for the complex conjugation on $\mathcal{A}$).

As explained in \cite{DH94}, first order differential calculi on a finite set are in bijective
correspondence with digraphs with at most a pair of edges connecting two vertices. Such a graph is characterized
by its adjacency matrix which is a square matrix $\mathcal{G}$ such $\mathcal{G}_{\m,\p}=1$ if there is a edge from $\m$ to $\p$ and $\mathcal{G}_{\m,\p}=0$ otherwise. Therefore, we can associate with each non-vanishing $\bb_{\m,\p}$ of some
differential calculus $( \Lambda^*\mathcal{A}, \dext)$ an
(undirected) edge connecting the nodes (vertices) $\m$ and $\p$ and hence the
adjacency matrix $\mathcal{G}$ is assigned by the sum of all 1-forms
\begin{eqnarray}
\label{G}{\bf G}=\sum_{\m \in \mathcal{L}}\sum_{j \in [n] \cup
[n]'}\bb_{\m,\m+\vv_j}=\sum_{j \in [n] \cup [n]'} \Theta^{\vv_j}.
\end{eqnarray}

The first order differential calculus then corresponds
to a complete graph where all vertices (the elements of $\mathcal{L}$) are connected to each other
by a pair of edges, whose interconnection structure is completely
encoded in the exterior derivative $\dext$. Indeed, from properties (\ref{simplex}) and
(\ref{transl}), it follows that our first order differential calculus is {\it inner} since $\G$ acts like $\dext$ when commuted with functions
$f \in \mathcal{A}$, i.e.

\begin{eqnarray}
\label{dbm}\left[ \G, f\right]&:=&\sum_{\m \in \mathcal{L}}f_m(\G \bb_{\m}-\bb_{\m}\G) \nonumber \\
&=&\sum_{\m \in \mathcal{L}}\sum_{j \in [n] \cup [n]'}f_m (\Theta^{\vv_j}\bb_{\m} -\bb_{\m}\Theta^{\vv_j})
\nonumber \\
&=&\sum_{\m \in \mathcal{L}}\sum_{j \in [n] \cup [n]'} f_m (\bb_{\m-\vv_j,\m} -\bb_{\m,\m+\vv_j}) \nonumber \\
&=&\dext f
\end{eqnarray}

To the {\it nearest-neighbour-nodes} being encoded in the 1-form
$\bb_{\m,\m+\vv_j}$, we associate a tangent vector field
$\partial_{\m+\vv_j,\m} \in \mathcal{T}_{\m}$.

From the $\dext-$action (\ref{diffun}) and the contraction
constraints (\ref{contraction}) and (\ref{contraction2}), we can
consider $\partial_{\m+\vv_j,\m} \in \mathcal{T}_{\m}$ as a partial
difference action at node $\m$:
\begin{eqnarray*}
\partial_{\m+\vv_j,\m}~f_{\m}=f_{\m+\vv_j}-f_{\m}, &
\mbox{for}~j=1,\ldots,n.
\end{eqnarray*}
Moreover, a basis $\{ \partial^{\vv_j} \}_{j \in [n]\cup [n]'}$ for
the vector space $\mathcal{T}$ satisfying the duality relation
\begin{eqnarray}
\label{dualityC}\langle \Theta^{\vv_k},\partial^{\vv_j}\rangle=\delta_{kj}
\end{eqnarray}
is uniquely determined by the relation
\begin{eqnarray}
\partial^{\vv_j}f=\sum_{\m \in \mathcal{L}}(\partial_{\m+\vv_j,\m}f_{\m})\bb_{\m}=T_{\vv_j}f-f.
\end{eqnarray}
Therefore,
\begin{eqnarray}
\label{df}\dext f=\left[ \G, f\right]= \sum_{j \in [n] \cup [n]'} (\partial^{\vv_j}f)\Theta^{\vv_j}
\end{eqnarray}

From the above characterization, we can now prove some results which relate $\dext$ with the graph
structure assigned by $\G$.

According to definition (\ref{prop4}) and property (\ref{nilp}), the set of $r-$forms is completely
determined by the set of edges of the graph. Indeed, the sum of all $r-$forms then corresponds to the $r-$power of $\G$,
\begin{eqnarray}
\label{Gr}
\G^r=\sum_{\{k_1,\ldots,k_r\} \subset [n] \cup
[n]'}\Theta^{\vv_{k_1}} \ldots \Theta^{\vv_{k_r}},
\end{eqnarray}
which is nothing else than the $r-$power of the adjacency matrix $\mathcal{G}$.

Let us take a close look at the sum of all $2-$forms. On the right-hand side of (\ref{Gr}), the summation over all indices $j,k \in [n] \cup [n]'$ involves terms of the form
\begin{eqnarray}
\label{Gr2}\frac{1}{2}\left(\Theta^{\vv_{j}} \Theta^{\vv_{k}}+\Theta^{\vv_{k}} \Theta^{\vv_{j}}\right)=\frac{1}{2}\sum_{\m \in \mathcal{L}} \left(\bb_{\m,\m+\vv_j,\m+\vv_j+\vv_k}+\bb_{\m,\m+\vv_k,\m+\vv_j+\vv_k} \right)
\end{eqnarray}
Because of $\bb_{\m,\m+\vv_j+\vv_k}=0$ the right-hand side of (\ref{Gr2}) in terms of the action of $\dext$ then corresponds to
\begin{eqnarray}
\label{Gr3}\frac{1}{2}\left(\Theta^{\vv_{j}} \Theta^{\vv_{k}}+\Theta^{\vv_{k}} \Theta^{\vv_{j}}\right)=\frac{1}{2}\sum_{\m \in \mathcal{L}} \dext \bb_{\m,\m+\vv_j+\vv_k}=0, &~j \neq k'.
\end{eqnarray}

Contrary to the oriented lattice reduction, in the symmetric lattice reduction the left-hand side of (\ref{Gr2}) does not vanish in general (see e.g. \cite{KK04}). However from (\ref{Gr3})
and
\begin{eqnarray*}
0=\sum_{\m \in \mathcal{L}}\dext \bb_{\m,\m}=-\sum_{j =
1}^n\left(\Theta^{\vv_{j}} \Theta^{\vv_{j'}}+\Theta^{\vv_{j'}} \Theta^{\vv_{j}}\right)
\end{eqnarray*}
the square of $\G$, $\G^2$, vanishes.

Therefore, to ensure that the term (\ref{Gr2}) vanishes for all $j,k \in [n] \cup [n]'$, our graphs should satisfy the following two conditions (see the paper by E. Forgy/U. Schreiber~\cite{FS04}):
\begin{enumerate}
\label{intedge}\item ${\bf G}$ has no {\it intermediate edges}:
$\bb_{\m,\p}\neq 0\Rightarrow \bb_{\m,\underline{l},\p}=0$,
\label{oppedge}\item ${\bf G}$ has no {\it opposite edges}:
$\bb_{\m,\p,\m}= 0.$
\end{enumerate}

Notice that on graphs without intermediate edges all the
2-paths of edges assigned by $\bb_{\m,\underline{l},\p}$ that connect the same two
points has to vanish. So our introduced lattice reduction
(\ref{latticeR}) is a particular case of a {\it graph without
intermediate edges}.

On the other hand, on graphs without {\it opposite edges}, neither multiple edges neither loops are admitted.

Therefore the above conditions leads the following theorem:
\begin{theorem}\label{extalgThm}
For graphs without {\it opposite edges}, the exterior
product rule
\begin{eqnarray*}
\label{extpalg}\left\{\Theta^{\vv_j},\Theta^{\vv_k}\right\}=0,
\end{eqnarray*}
holds for all $j,k \in [n]\cup [n]'$.
\end{theorem}
Here and elsewhere, $\{ a,b\}=ab+ba$ denotes the anti-commutator bracket of $a$ and $b$.

By applying the coordinate definition (\ref{df}), Theorem (\ref{extalgThm}) and the nil-potency of $\dext$ (\ref{nilp}), we arrive at the
following corollary:
\begin{corollary} \label{commuteT}
For graphs without {\it opposite edges}, the basis elements
$\partial^{\vv_j},j \in [n] \cup [n]'$, mutually commute when acting
on functions $\mathcal{A}$, that is
\begin{eqnarray*}
\partial^{\vv_j}(\partial^{\vv_k}f)=\partial^{\vv_k}(\partial^{\vv_j}f), & \mbox{for all}~j,k \in [n]\cup [n]'~\mbox{and}~f\in \mathcal{A}.
\end{eqnarray*}
\end{corollary}

Let us remark that, for the product of two functions $f,g \in
\mathcal{A}$, the difference action $\partial^{\vv j}$ satisfies the
product rule
\begin{eqnarray}
\label{leibniz2}\partial^{\vv_j}(f
g)=(\partial^{\vv_j}f)(T_{\vv_j}g)+f(\partial^{\vv_j}g).
\end{eqnarray}

This means that the finite difference actions $\partial^{\vv_j}$ do
not obey the ordinary Leibniz rule. In fact, the application of
$\partial^{\vv_j}$ to, say, higher powers becomes increasingly
cumbersome.

Due to the discreteness of the formalism and, as a consequence, the
inevitable bi-locality of $\partial^{\vv_j}$ there is no chance to
get something as a `true' Leibniz rule on this level. Nevertheless,
there is a certain systematic in it, namely the product rule.

We will come back to the non-Leibniz character of
$\partial^{\vv_j}$, when establishing the duality between
differential forms and vector-fields. Before however doing that we
will need to further clarify the role of $\dext$ on $r-$forms.

Under the conditions of Theorem~\ref{extalgThm} we get
$$
\bb_{\m}\G \Theta^{\vv_j}=-\bb_{\m,\m+\vv_j}\G=-\sum_{k \in [n] \cup [n]'}
\bb_{\m,\m+\vv_j,\m+\vv_j+\vv_k}.$$ Then we have
$$\dext \bb_{\m,\m+\vv_j}=[\G,\bb_{\m}]\Theta^{\vv_j}=\G \bb_{\m,\m+\vv_j}+\bb_{\m,\m+\vv_j}\G.$$
Applying induction on $r>0$ and the graded Leibniz rule $$\dext
(\omega_1 \omega_r)=(\dext \omega_1)\omega_r-\omega_1(\dext
\omega_r)$$ together with the action of $\dext$ on $r$-forms
(\ref{prop4}), lead to the following theorem.
\begin{theorem}\label{dwr}
For graphs without {\it opposite edges}, the actions of $\G$ and $\dext$ on $r-$forms are related as follows:
\begin{eqnarray*}
\dext \omega_r &=&\G\omega_r-(-1)^r\omega_r \G
\end{eqnarray*}
\end{theorem}

\subsection{Differential forms representations of $End(\Lambda^*\mathcal{A})$ on the symmetric lattice}

In order to introduce directions in lattices with mesh-width $h>0$,
we introduce the coordinate functions
\begin{eqnarray}
\label{xj}x_j&=&\sum_{\m \in \mathcal{L}}h~m_j\bb_{\m}
\end{eqnarray}
and split up $\G$ into $\G=\overleftarrow{\G}+\overrightarrow{\G}$,
where
\begin{eqnarray*}
\label{Gfb}\begin{array}{ccc}
\label{Gbackward}\overleftarrow{\G}=\sum_{j =1}^n\Theta^{\vv_{j}'}, &
\label{Gforward}\overrightarrow{\G}=\sum_{j=1}^n \Theta^{\vv_{j}}.
\end{array}
\end{eqnarray*}
Hence, by applying $\dext$ on $x_j$, we can find the coordinate
differentials $\dext x_j$ on the symmetric lattice as
\begin{eqnarray}
\label{dxj}\dext x_j =[\G,x_j]=-\dext x_j^-+\dext x_j^+ .
\end{eqnarray}
where $\dext x_j^\pm$ corresponds to
\begin{eqnarray*}
\label{dxj2}\dext x_j^+ = \left[\overrightarrow{\G},x_j\right]=h\Theta^{\vv_{j}}, &
 \dext x_j^-=\left[-\overleftarrow{\G},x_j\right]=h\Theta^{\vv_{j'}}.
\end{eqnarray*}

According to (\ref{xj}) we can express any function $f \in
\mathcal{A}$ as $f(\x)$ and by direct application of the coordinate
differentials (\ref{dxj2}), the non-commutativity of functions and
$1-$forms shown in (\ref{transl}) can now be given in terms of
$\underline{x}=(x_1,\ldots,x_n)$ as
\begin{eqnarray}
\label{transl2}\dext x_j^\pm f(\x)=f({\x\pm h\vv_{j}})\dext x_j^\pm, &j = 1,\ldots,n.
\end{eqnarray}

To represent any r-form $\omega_r \in \Lambda^*\mathcal{A}$ in terms
of the coordinate differentials $\dext x_j^{\pm}$, we proceed as
follows:

For two ordered subsets $J'\subset[n]'$ and $K \subset[n]$ given by
$$J'=\{j_1',\ldots,j_s'\},~~K=\{k_{1},\ldots, k_{t} \}\subset [n] ,$$ we
put
\begin{eqnarray*}
\begin{array}{ccc}
\dext x_{\emptyset}^\pm=1, \\
\dext x_J^-=\dext x_{j_1}^- \ldots \dext x_{j_s}^- \\
\dext x_K^+=\dext x_{k_1}^+ \ldots \dext x_{k_t}^+.
\end{array}
\end{eqnarray*}
Therefore, the $r-$form defined by (\ref{rform}) then corresponds to
\begin{eqnarray*}
\label{wr}\omega_r(\x)&=&\sum_{|J'|+|K| =r} F_{J'K}(\x)\dext
x_J^-\dext x_K^+,
\end{eqnarray*}
where $F_{J'K}$ denotes the function (\ref{function}) indexed by $J' \cup K$ and $\sum_{|J'|+|K| =r}$ denotes a sum restricted to ordered subsets $J'\subset[n]'$ and $K \subset[n]$.

For the vector fields $\partial_h^{\pm j}\in \mathcal{T}$ defined as the forward/backward difference actions
\begin{eqnarray}
\label{finitediff}(\partial_h^{+j}f)(\x)=
\frac{1}{h}(\partial^{\vv_j}f)(\x), & (\partial_h^{-j}f)(\x)=
-\frac{1}{h}(\partial^{\vv_{j'}}f)(\x), & \forall {f \in
\mathcal{A}}
\end{eqnarray}
we can split the $\dext-$action on $\omega_r$ as
\begin{eqnarray}
\label{dG}\dext \omega_r(\x)&=&  \dext_+\omega_r(\x)-\dext_-\omega_r(\x),
\end{eqnarray}
where $\dext_\pm \omega_r(\x)$ are given by
\begin{eqnarray}
\label{dpm}\dext_{\pm} \omega_r(\x)&=& \sum_{j=1}^n\sum_{|J'|+|K|
=r}(\partial_h^{\pm j}F_{J'K})(\x)~\dext x_j^{\pm}\dext x_J^-\dext
x_K^+.
\end{eqnarray}

The above framework suggests the following bi-graded algebra (i.e. a
{\it bi-complex}) decomposition
\begin{eqnarray*}
\Lambda^*\mathcal{A}= \sum_{p,q=0}^n\bigoplus \Lambda^{p,q}\mathcal{A}
\end{eqnarray*}
where $\mathcal{A}:=\Lambda^{0,0}\mathcal{A}$ and the discrete
exterior differential maps
\begin{eqnarray*}
\dext_-: \Lambda^{p,q}\mathcal{A} \rightarrow \Lambda^{p+1,q}\mathcal{A}, & \dext_+: \Lambda^{p,q}\mathcal{A} \rightarrow \Lambda^{p,q+1}\mathcal{A}
\end{eqnarray*}
are defined by formula (\ref{dpm}). Under the conditions of Theorem~\ref{extalgThm}, we
have that the coordinate differentials $\dext x_j^\pm$ satisfy the
anti-commutation relations
\begin{eqnarray}
\label{grassman}\begin{array}{ccc}\left\{\dext x_j^{\pm},\dext x_k^{\pm}\right\} & =0
& \forall_{j,k =1,\ldots,n}, \\
\left\{\dext x_j^{+},\dext x_k^{-}\right\} & =0
& \forall_{j,k =1,\ldots,n}.
\end{array}
\end{eqnarray}
and, furthermore,
\begin{eqnarray*}
\label{nilp2}\dext_{\pm}(\dext_{\pm}\omega_r(\x))=0= \dext_+(\dext_-\omega_r(\x))+\dext_-(\dext_+\omega_r(\x)).
\end{eqnarray*}

The basic endomorphisms acting on $\Lambda^*\mathcal{A}$ in an
exterior way are the linear operators $\gamma^{\pm j}\in
\mbox{End}(\Lambda^*\mathcal{A})$ defined as
\begin{eqnarray}
\label{exterior}\begin{array}{ccc}\gamma^{-j}:\Lambda^{p,q}\mathcal{A} \rightarrow \Lambda^{p+1,q}\mathcal{A}, &\omega(\x) \mapsto \dext x_j^{-}\omega(\x), \\
\gamma^{+j}:\Lambda^{p,q}\mathcal{A} \rightarrow \Lambda^{p,q+1}\mathcal{A}, &\omega(\x) \mapsto \dext x_j^{+}\omega(\x).
\end{array}
\end{eqnarray}
Having defined the left and right exterior product representations,
it raises the question how to define left and right representations
for the interior products $\vartheta^{\pm j} \in
\mbox{End}(\Lambda^*\mathcal{A})$ in terms of the duality
contraction $\langle \cdot,\cdot \rangle$. To be consistent with the
nature of the interior product, we impose the duality conditions
$\vartheta^{\mp j}(\dext x_k^{\pm})=0$ and $\vartheta^{\pm j}(\dext
x_k^{\pm})=\delta_{jk}$.

Let us proceed as follows: first since the Leibniz rule
(\ref{leibniz}) is also valid for the exterior derivatives
$\dext_\pm$ we observe
\begin{eqnarray*}
\dext_{\pm}(x_k \omega_r(\x))=\dext x_k^{\pm}\omega_r(\x)+x_k\dext_{\pm}\omega_r(\x)
\end{eqnarray*}
and, hence, direct application of the duality contraction property
(\ref{dualityC}) and the non-commutativity (\ref{transl2}) lead to
\begin{eqnarray}
\label{productC}\begin{array}{ccc}\left\langle\dext x_k^{\pm}\omega(\x),\partial_h^{\pm j}\right\rangle&=& \delta_{jk}(T_h^{\pm j}\omega)(\x)-\dext x_k^{\pm}\left\langle\omega(\x),\partial_h^{\pm j}\right\rangle, \\
\left\langle\dext x_k^{\mp}\omega(\x),\partial_h^{\pm
j}\right\rangle&=& -\dext
x_k^{\mp}\left\langle\omega(\x),\partial_h^{\pm j}\right\rangle.
\end{array}
\end{eqnarray}
Hereby, $T_h^{\pm j}$ denotes the shift action $(T_h^{\pm
j}\omega)(\x)=\omega(\x\pm h \vv_j)$ on the exterior algebra
$\Lambda^*\mathcal{A}$.

This suggests that the interior product operators $\vartheta^{\pm j}
\in \mbox{End}(\Lambda^*\mathcal{A})$ should be defined as
contraction operators with shifting role opposite to the
differential form
\begin{eqnarray}
\label{interiorP}\begin{array}{ccc}\vartheta^{-j}:\Lambda^{p,q}\mathcal{A} \rightarrow \Lambda^{p-1,q}\mathcal{A}, &\omega(\x) \mapsto \left\langle (T_h^{+j}\omega)(\x),\partial_h^{-j}\right\rangle, \\
\vartheta^{+j}:\Lambda^{p,q}\mathcal{A} \rightarrow \Lambda^{p,q-1}\mathcal{A}, &\omega(\x) \mapsto \left\langle (T_h^{-j}\omega)(\x),\partial_h^{+j}\right\rangle.
\end{array}
\end{eqnarray}
There are some formulae that follow from
equations (\ref{grassman}),~(\ref{nilp2}) and (\ref{interiorP})
will be of interest, namely, the anti-commutation relations between
the interior and the exterior product representations $\gamma^{\pm
j},\vartheta^{\pm j}\in \mbox{End}(\Lambda^*\mathcal{A})$
\begin{eqnarray}
\label{fermionic}\begin{array}{ccc}
\gamma^{\pm j}(\gamma^{\pm k}\omega(\x))+\gamma^{\pm k}(\gamma^{\pm j}\omega(\x))&=&0, \\
\gamma^{+j}(\gamma^{- k}\omega(\x))+\gamma^{-k}(\gamma^{+ j}\omega(\x))&=&0, \\
\vartheta^{\pm j}(\vartheta^{\pm k}\omega(\x))+\vartheta^{\pm k}(\vartheta^{\pm j}\omega(\x))&=&0, \\
\vartheta^{+j}(\vartheta^{- k}\omega(\x))+\vartheta^{-k}(\vartheta^{+ j}\omega(\x))&=&0, \\
\gamma^{+j}(\vartheta^{- k}\omega(\x))+\vartheta^{-k}(\gamma^{+ j}\omega(\x))&=&0, \\
\gamma^{\pm j}(\vartheta^{\pm k}\omega(\x))+\vartheta^{\pm
k}(\gamma^{\pm j}\omega(\x))&=&\delta_{jk}\omega(\x). \end{array}
\end{eqnarray}

The above identities then correspond to the graded algebra like the one of Fermionic creation and annihilation operators \cite{nielsen91}.




It's now interesting to compare our approach with the approach proposed by E.~Forgy/U.~Schreiber in \cite{FS04}, where they construct a discrete differential calculi on causal graph complexes.

On their terminology (see \cite{FS04}, page 4), the algebra of endomorphisms $\mbox{End}(\Lambda^*\mathcal{A})$ would be the inner product space $\mathcal{H}(\mathcal{A}, \langle \cdot |\cdot \rangle)$. Indeed a representation of the basic endomorphisms $\vartheta^{\pm j}$ in terms of the inner product $\langle \cdot | \cdot \rangle$ already tacitly exists (see \cite{FS04}, section 3.6). However, bear in mind that these operators can be introduced independently of the inner product $\langle \cdot |\cdot \rangle$ and the behavior is invariant under diffeomorphisms (c.f. \cite{rota85}).

Moreover, our framework gives a quite natural way to describe the algebra $\mbox{End}(\Lambda^*\mathcal{A})$ as a canonical realization of the
Fermi algebra in a metrically independent, via duality arguments.

We are now in conditions to explore the correspondence between
differential forms and Clifford algebras on the symmetric lattice.

\subsection{Clifford algebras and hermitian structure of the symmetric lattice}\label{CliffordLattice}

As it was shown in~\cite{BDS82,GM92,Sommen97} Clifford algebras can
be defined in several different ways. One of these ways is as a
subalgebra of the algebra of endomorphisms of the exterior algebra.

In order to get this correspondence for the symmetric lattice, we
proceed as follows:

First of all, notice that from the commutation relations
(\ref{fermionic}) the endomorphisms
$\xi^{+j},\xi^{-j}~:~j=1,\ldots,n,\in
\mbox{End}(\Lambda^*\mathcal{A})$ defined by
\begin{eqnarray}
\label{xipm}\xi^{\pm j}=\gamma^{\pm j}+\vartheta^{\mp j}
\end{eqnarray}
satisfy the graded fermionic identities when acting on
$\Lambda^*\mathcal{A}:$
\begin{eqnarray*}
\label{wittdiff}\begin{array}{ccc}\xi^{\pm j}(\xi^{\pm k}\omega(\x))+\xi^{\pm k}(\xi^{\pm
j}\omega(\x))&=&0, \\
\xi^{+j}(\xi^{-k}\omega(\x))+\xi^{-k}(\xi^{+j}\omega(\x))&=&\delta_{jk}\omega(\x)
\end{array} & \mbox{for all}~{\omega \in \Lambda^*\mathcal{A}}.
\end{eqnarray*}
Furthermore, the elements
\begin{eqnarray}
\label{Upm}\Upsilon^{\pm j}=\xi^{+j}\pm \xi^{-j}
\end{eqnarray}
satisfy the graded orthogonal identities when acting on
$\Lambda^*\mathcal{A}:$
\begin{eqnarray*}
\begin{array}{ccc}\Upsilon^{\pm j}(\Upsilon^{\pm k}\omega(\x))+\Upsilon^{\pm k}(\Upsilon^{\pm
j}\omega(\x))&=&\pm 2\delta_{jk}\omega(\x), \\
\Upsilon^{+j}(\Upsilon^{-k}\omega(\x))+\Upsilon^{-k}(\Upsilon^{+j}\omega(\x))&=&0
\end{array} & \mbox{for all}~{\omega \in \Lambda^*\mathcal{A}}.
\end{eqnarray*}
This clearly suggests that $\Upsilon^{\pm j}$ behave like the
generators of the real Clifford algebra of signature $(n,n)$,
$\BR_{n,n}$.

For the symmetric lattice we get the following interesting features
of the coordinate differentials $\dext x_j=\Upsilon^{-j}({\bf 1})$
and $\dext \tau_j:=\dext x_j^++\dext x_j^-=\Upsilon^{+j}({\bf 1})$:
\begin{eqnarray*}
(\dext x_j)'=-\dext x_j, & (\dext \tau_j)'=\dext \tau_j
\\
(\dext x_j)^{\widetilde{~}}=\dext x_j, & (\dext
\tau_j)^{\widetilde{~}}=-\dext \tau_j \\
 (\dext x_j)^{\dag}=-\dext x_j, &
(\dext \tau_j)^{\dag}=\dext \tau_j.
\end{eqnarray*}
Hence $\dext x_j$ behaves as a real while $\dext \tau_j$ behaves as pure imaginary
with respect to the involution, reversion and $\dag-$conjugation, respectively. It then turns out that the automorphisms
(\ref{involution})-(\ref{hconjugation}) play the same role as the
automorphisms on the complex Clifford algebra $\BC_{2n}$ according
to the symmetric nature of the lattice.

Since the real Clifford algebra $\BR_{n,n}$ is contained in the
complex Clifford algebra $\BC_{2n}=\BC \bigotimes \BR_{0,2n}$
as a special subalgebra \cite{BDS82}, the isomorphism between the
$\mbox{End}(\Lambda^*\mathcal{A})$ and $\BC_{2n}=\BC \bigotimes
\BR_{0,2n}$~is thus obtained through the identification $\e_j
\leftrightarrow \Upsilon^{-j}$ and $\e_{j+n} \leftrightarrow
i\Upsilon^{+j}$.

The corresponding Witt basis for $\BC_{2n}$ is given by
\begin{eqnarray*}
\f_j=\frac{1}{2}(\e_j-i\e_{n+j}), &
\f_j^\dag=-\frac{1}{2}(\e_j+i\e_{n+j})
\end{eqnarray*}
and satisfies the following anti-commuting identities
\begin{eqnarray*}
\mbox{Grassmann identities:} & \{\f_j,\f_k\}=0=\{\f_j^\dag,\f_k^\dag\}, &  \\
\mbox{duality identities:} &\{\f_j,\f_k^\dag\}=\delta_{jk}.
\end{eqnarray*}
Moreover, we can identify it with the set of endomorphisms $$\{
\xi^{+j},\xi^{-j}~:~j=1,\ldots,n\}\subset
\mbox{End}(\Lambda^*\mathcal{A})$$ {\it viz} $\f_j \leftrightarrow
\xi^{+j}$ and $\f_j^\dag \leftrightarrow \xi^{-j}$.

All the above identifications clearly suggests the hermitian
Clifford setting as the natural multi-vector setting to develop a
discrete framework on the symmetric lattice.

It is now interesting to compare our approach with the approach proposed by Kanamori and Kawamoto in~\cite{KK04}.

We would like to remark that the initial construction proposed in
subsection \ref{graphS} is analogous to the approach proposed in
that paper. However, while the authors define a Clifford-like
algebra on the symmetric lattice by introducing a new Clifford
product in a rather ``artificial'' way, we avoid this problem by
passing to the algebra of endomorphisms.

While in the approach of Kanamori and Kawamoto in~\cite{KK04} there
is no restriction to obtain Clifford products which are associative
and distributive, however, that construction is metric-dependent,
since the Clifford product can be indefinite.

In our approach, we note that {\it a priori} no problem in relation
to the question of associativity and distributivity on the lattice
will appear since the algebra of endomorphisms of a given space
equipped with the standard sum and a product defined by composition
is obviously associative and distributive. On the other hand, the
structure of the algebra of endomorphisms lead also to a
construction of Clifford products which are metric-independent in
the sense that they only depend on the duality between the tangent
and the cotangent space.

\section{Dirac operators and vector variables on the symmetric lattice}

The main objective in this section is to show some similarities
between the hermitian setting and the symmetric structure of the lattice.

Our starting point again is the definition of $\dext_\pm$ and $\dext=\dext_+ -\dext_-$.

From Corollary~\ref{commuteT} we know that all forward and backward
differences mutually commute when acting on $\Lambda^*\mathcal{A}$
\begin{eqnarray*}
\label{annihilation}\begin{array}{ccc}\partial_h^{\pm j}(\partial_h^{\pm k}\omega)(\x)=\partial_h^{\pm k}(\partial_h^{\pm j}\omega)(\x)\\
 \partial_h^{+ j}(\partial_h^{- k}\omega)(\x)=\partial_h^{- k}(\partial_h^{+ j}\omega)(\x),
 \end{array}& \forall {\omega \in  \Lambda^*\mathcal{A}}.
\end{eqnarray*}
Furthermore, they are interrelated by the translations $(T_h^{\pm
j}\omega)(\x)=\omega(\x\pm h \vv_j)$
\begin{eqnarray*}
T_h^{-j}(\partial_h^{+j}\omega)(\x)=(\partial_h^{-j}\omega)(\x), & T_h^{+j}(\partial_h^{-j}\omega)(\x)=(\partial_h^{+j}\omega)(\x).
\end{eqnarray*}

Using (\ref{transl2}) and (\ref{dpm}), we can thus write the
exterior differentials $\dext_{\pm}$ in the form
\begin{eqnarray*}
\dext_{\pm}=\sum_{j=1}^n \dext x_j^{\pm}\partial_h^{\mp j}.
\end{eqnarray*}
Let us now introduce the symmetric and skew-symmetric difference
operators $\nabla_h^j$ and ${\tilde\nabla}_h^j$, respectively, as:
\begin{eqnarray*}
\label{diffops}\nabla^{j}_h=\frac{1}{2}(\partial_h^{-j}+\partial_h^{+j})
,&
{\tilde\nabla}^{j}_h=\frac{1}{2i}(\partial_h^{-j}-\partial_h^{+j}).
\end{eqnarray*}
Then the operator $\dext$ defined in (\ref{dG}) corresponds in
terms of the coordinate differentials $\dext x_j$ and $\dext \tau_j$
to
\begin{eqnarray*}
\dext&=&\dext_+-\dext_-\nonumber \\
&=&\sum_{j=1}^n \dext x_j^+\partial_h^{-j}-\dext x_j^-\partial_h^{+j} \\
&=&\sum_{j=1}^n \dext x_j \nabla^{j}_h+i \dext \tau_j
{\tilde\nabla}^{j}_h.
\end{eqnarray*}
The above identities suggest the introduction of the following
operators acting on $End(\Lambda^*\mathcal{A})$
\begin{eqnarray}
\label{Dpm}\partial_\pm&=&\sum_{j=1}^n \xi^{\pm j} \partial_h^{\mp j} \\
\label{D}\partial&=&\sum_{j=1}^n \Upsilon^{-j}\nabla^{j}_h+ i\Upsilon^{+j}{\tilde\nabla}^{j}_h
\end{eqnarray}
where $\xi^{\pm j}$ and $\Upsilon^{\pm j}$ are the basic endomorphisms defined in (\ref{xipm}) and (\ref{Upm}), respectively.
All the geometry of the symmetric lattice is now encoded in the operators $\partial_\pm$ and $\partial$, which are nothing else than the
 hermitian and the orthogonal Dirac-K\"ahler operators acting on symmetric lattices, respectively.

 The operator defined in (\ref{D}) is closely related the Dirac operator introduced by Kanamori, Kanamoto
 in \cite{KK04} (see formula (4.10), page 21) and Forgy, Schreiber in \cite{FS04} (see formula (5.34), page 71).

Here we are interested to establish a correspondence with the Hermitian setting proposed in \cite{BSS07}, our departure point will be the operators (\ref{Dpm})
instead of the operator (\ref{D}). This is indeed the main difference between our approach and the approaches proposed by Kanamori, Kanamoto in
\cite{KK04} and Forgy, Schreiber in \cite{FS04}.

From the correspondence $\f_j \leftrightarrow
\xi^{+j}$ and $\f_j^\dag \leftrightarrow \xi^{-j}$, we see that $\partial_{\pm}$ is the lattice counterpart to the hermitian Dirac operator and its conjugate on $\BC_{2n}$
\begin{eqnarray}
\label{hDirac}\begin{array}{ccc}\partial_{z}=\sum_{j=1}^n \f_j \partial_h^{-j}, \\
\partial_{z}^\dag=\sum_{j=1}^n \f_j^\dag \partial_h^{+j}
\end{array}
\end{eqnarray}
where $\partial_h^{\pm j}=\frac{1}{2}\left(\nabla_h^j \mp i
{\tilde\nabla}_h^{j}\right)$ are the discrete counterparts of the
classical Cauchy-Riemann operators and their conjugates. Moreover,
from the correspondence $\e_j \leftrightarrow \Upsilon^{-j}$,
$\e_{j+n} \leftrightarrow i\Upsilon^{+j}$ we see that $\partial$ is
the lattice counterpart of the discrete Dirac operator on
$\BC_{2n}$
\begin{eqnarray}
\label{Dirac}\partial_{X}&=&\sum_{j=1}^n \e_j\nabla^{j}_h+
\e_{j+n}{\tilde\nabla}^{j}_h
\end{eqnarray}
Similarly, the lattice counterpart $-i(\partial_{+}+\partial_{-})$ of the discrete Dirac operator on $\BC_{2n}$ is given by
\begin{eqnarray}
\label{Dirac2}\partial_{X|}:=-i(\partial_{z}+\partial_{z}^\dag)&=&\sum_{j=1}^n
\e_j{\tilde\nabla}^{j}_h- \e_{j+n}\nabla^{j}_h
\end{eqnarray}
There are some formulae holding for the operators (\ref{hDirac}),(\ref{Dirac}) and (\ref{Dirac2}) that
will be of interest in the hermitian setting, namely
\begin{enumerate}
\item Isotropy condition:  $ \partial_z^{2}=0=(\partial_z^{\dag})^2$
\item Orthogonality condition: $\{\partial_X,\partial_{X|}\}=0$
\item Star Laplacian splitting:
\begin{itemize}
\item Using Dirac operators: $\partial_X^2=-\sum_{j=1}^n\partial_h^{-j}\partial_h^{+j}=\partial_{X|}^2$
\item Using hermitian Dirac operators: $\{\partial_z,\partial_z^\dag\}=\sum_{j=1}^n\partial_h^{-j}\partial_h^{+j}$
\end{itemize}
\end{enumerate}

Having established the correspondence between hermitian Dirac and Dirac
operators on the symmetric lattice and discrete Dirac operators on
the Clifford algebra, we arrive now at the question how to make the
correspondence between coordinate vector functions on the symmetric
lattice and vector variables on the Clifford algebra.

Let us remark that the finite difference action $\partial_h^{\pm j}$
acting on $\Lambda^*\mathcal{A}$ satisfies the product rule
\begin{eqnarray}
\label{duality}\partial_h^{\pm j}\left(x_k  (T_h^{\mp k}\omega)(\x)\right)&=&\delta_{jk}T_h^{\mp k}(T_h^{\pm j}\omega)(\x)+x_kT_h^{\mp k}(\partial_h^{\pm j}\omega)(\x) \nonumber \\
&=&\delta_{jk}\omega(\x)+x_kT_h^{\mp k}(\partial_h^{\pm j}\omega)(\x).
\end{eqnarray}
which establishes the duality between the finite difference
operators $\partial_{h}^{\pm j}$ and the ``formal'' coordinate
functions $x_jT_{h}^{\mp j}$ on $\mbox{End}(\Lambda^*\mathcal{A})$.

We also note that the coordinate variables $x_j T_h^{\pm j}$ mutually commute, when acting on functions on $ \Lambda^*\mathcal{A}$
\begin{eqnarray}
\label{creation}x_jT_h^{\pm j}(x_kT_h^{\pm k}\omega(\x))=x_kT_h^{\pm
k}(x_jT_h^{\pm j}\omega(\x)). \end{eqnarray} Furthermore, the commutative
relations (\ref{annihilation}),(\ref{creation}) together with the
duality relation (\ref{duality}) endow an algebraic representation
of the Weyl-Heisenberg algebra~\cite{nielsen91}, where the
``formal'' coordinate functions $x_jT_h^{\pm j}$ represent
``creation'' operators dual to the ``annihilation'' operators
$\partial_h^{\pm j}$.

From the previous relations, (\ref{productC}), and definition
(\ref{interiorP}), we immediately obtain the following commutation
relations between $\partial^{\pm j}_h$ and $\gamma^{\pm j}$ when
acting on $\mbox{End}(\Lambda^*\mathcal{A})$
\begin{eqnarray*}
\begin{array}{ccc}
\partial^{\pm j}_h(\gamma^{\pm k}(\omega(x)))=\gamma^{\pm k}(\partial^{\pm j}_h(\omega(x))), \\
\partial^{\pm j}_h(\gamma^{\mp k}(\omega(x)))=\gamma^{\mp k}(\partial^{\pm
j}_h(\omega(x))).
\end{array}
\end{eqnarray*}
Moreover, Clifford-like operators on the symmetric lattice are
therefore encoded in the following algebra of endomorphisms
\begin{eqnarray*}
Alg\left\{\partial_h^{-j},\partial_h^{+j},x_jT_h^{+j},x_jT_h^{-j},\xi^{+j},\xi^{-j}~:~j=1,\ldots,n\right\},
\end{eqnarray*}
where $\xi^{\pm j}\in \mbox{End}(\Lambda^*\mathcal{A})$ defined in (\ref{xipm}) satisfy the fermionic relations (\ref{wittdiff}).

Let us define formally $z_j:=x_jT_h^{+j},~ \overline{z_j}:= x_jT_h^{-j}$
as the complex variables $z_j$ and their conjugates
$\overline{z_j}$, respectively. Defining the hermitian vector variable $z=\sum_{j=1}^n \f_j z_j$
and its hermitian conjugate $z^\dag=\sum_{j=1}^n \f_j^\dag
\overline{z_j}$, the Clifford vector variable $X$ associated to
$\partial_{X}$ takes the form $X=z-z^\dag$ while the Clifford vector
variable $X|$ associated to $\partial_{X|}$ takes the form
$X|=-i(z+z^\dag)$.

Hence, the following ``formal'' vector variable identifications
naturally follows:
\begin{eqnarray*}\label{vectorV}\begin{array}{ccc}
z \longleftrightarrow \sum_{j=1}^n \xi^{+j} (x_jT^{+j}_h) , \qquad
z^\dag\longleftrightarrow \sum_{j=1}^n \xi^{-j} (x_jT^{-j}_h), \\ \ \\
X \longleftrightarrow \sum_{j=1}^n \Upsilon^{-j}~ \frac{x_j(T^{+j}_h+T^{-j}_h)}{2}+i\Upsilon^{+j}~ \frac{x_j(T^{+j}_h-T^{-j}_h)}{2i}, \\ \ \\
X| \longleftrightarrow \sum_{j=1}^n \Upsilon^{-j}~
\frac{x_j(T^{+j}_h-T^{-j}_h)}{2i}-i\Upsilon^{+j}\frac{x_j(T^{+j}_h+T^{-j}_h)}{2}~.
\end{array}
\end{eqnarray*}
There are some formulae holding for the vector variables
$X$,$X|$, $z$, and $z^\dag$ that will be of interest in the hermitian setting, namely
\begin{enumerate}
\item The isotropy condition:  $ z^2=0=(z^{\dag})^2$
\item The orthogonality condition: $\{X,X|\}=0$
\item Square variable splitting:
\begin{itemize}
\item Using Clifford vector variables: $$X^2=-\sum_{j=1}^n (x_jT_h^{+j})(x_jT_h^{-j})=(X|)^2$$
\item Using hermitian vector variables: $$\{z,z^\dag\}=\sum_{j=1}^n (x_jT_h^{+j})(x_jT_h^{-j})$$
\end{itemize}
\end{enumerate}
By Clifford geometric product \cite{BDS82,BDS05},
$$ ab=a \bullet b+a \wedge b$$ one can also introduce the discrete hermitian Euler operators as $E_{z}=2 z \bullet \partial_z$,~$E_{z^\dag}=2 z^\dag \bullet \partial_z^\dag$ and the hermitian Gamma operators as $\Gamma_z=z \wedge \partial_z$, $\Gamma_{z^\dag}=z^\dag \wedge \partial_z^\dag$.

The Weyl-Heisenberg character of the operators $\partial_h^{\pm j}$ and $x_jT^{\mp j}_h$ then lead to the following elegant formulae similar to the intertwining relations concerning the operator $\partial_z,\partial_z^\dag,z,z^\dag,E_{z},E_{z^\dag},\Gamma_{z}$ and $\Gamma_{z^\dag}$ in \cite{Bures07,BSS07}.
\begin{eqnarray}
\label{hformulae}\begin{array}{ccc}
\{ z,\partial_z\}=\beta+E_z, & [ z,\partial_z]=-\beta+\Gamma_z
\\
\{ z^\dag,\partial_z^\dag\}=(n-\beta)+E_{z^\dag}, & [z^\dag,\partial_z^\dag]=-(n-\beta)+\Gamma_{z^\dag}\\
\{ z^\dag,\partial_z\}=0, & \{z,\partial_z^\dag\}=0
\end{array}
\end{eqnarray}
where $\beta=\sum_{j=1}^n \f_j^\dag\f_j$ denotes the spin-Euler operator.

Analogously, the link between the hermitian operators $E_{z},E_{z^\dag},~\Gamma_{z},~\Gamma_{z^\dag}$ and the traditional Euler and Gamma operators $E_X=-X \bullet \partial_{X},E_{X|}=-X| \bullet \partial_{X|},~\Gamma_{X}=-X \wedge \partial_{X}$ and $\Gamma_{X|}=-X| \wedge \partial_{X|}$ can also be obtained:
\begin{eqnarray*}
\begin{array}{ccc}E_X=E_{z}+E_{z^\dag}=E_{X|}, \\
\Gamma_X=\Gamma_z+\Gamma_{z^\dag}-2( z^\dag \wedge \partial_{z}+ z \wedge \partial_{z}^\dag), \\
\Gamma_X=\Gamma_z+\Gamma_{z^\dag}+2( z^\dag \wedge \partial_{z}+ z \wedge \partial_{z}^\dag).
\end{array}
\end{eqnarray*}
Let us take now a close look at the concept of discrete homogeneity
on $\BC_{2n}$.

As homogeneous polynomials are expected to be $\BC_{2n}-$valued
eigenfunctions of the Euler operator corresponding to the eigenvalue given by the degree
of the polynomial, in the hermitian setting it still makes sense to define
discrete homogeneous polynomials of degree $(p,q)$ as solutions of
the coupled eigenvalue problem
\begin{eqnarray}
\label{coupleEig}\left\{\begin{array}{ccc}
E_z [R_{p,q}(z,z^\dag)]=pR_{p,q}(z,z^\dag) \\ \ \\
 E_{z^\dag} [R_{p,q}(z,z^\dag)]=qR_{p,q}(z,z^\dag)
 \end{array}\right.
\end{eqnarray}
and discrete hermitian monogenic homogeneous polynomials of degree $(p,q)$
as the solutions of (\ref{coupleEig}) satisfying the zero hermitian Dirac
constraints:
\begin{eqnarray}
\label{hmon}\partial_zR_{p,q}(z,z^\dag)=0=\partial_{z^\dag}R_{p,q}(z,z^\dag)
\end{eqnarray}
A further consequence of (\ref{hformulae})-(\ref{hmon}), is that
any discrete hermitian monogenic homogeneous polynomial of degree $(p,q)$ is
an eigenfunction of the hermitian Gamma operators, namely $\Gamma_z
[R_{p,q}(z,z^\dag)]=-pR_{p,q}(z,z^\dag)$ and $\Gamma_{z^\dag}
[R_{p,q}(z,z^\dag)]=-qR_{p,q}(z,z^\dag)$.

It is interesting to see the similarities between the solutions of
the coupled eigenvalue problem (\ref{coupleEig}) and the concept of
discrete homogeneous polynomials introduced in our previous
paper~\cite{FK07}.

In the terminology of that paper, $E_h^{\pm}=\sum_{j=1}^n x_j
\partial_h^{\pm j}$ is the forward/backward Euler operator while
$(x)_{\pm}^{(\alpha)}$ are the multi-index factorial powers
satisfying the eigenvalue property
\begin{eqnarray*}
E_h^{\pm}(x)_{\pm}^{(\alpha)}=|\alpha|(x)_{\pm}^{(\alpha)}.
\end{eqnarray*}
On the other hand, discrete hermitian Euler operators coincide with the forward/backward difference Euler operators, namely $E_{z}=E_h^+$ and $E_{z^\dag}=E_h^-$.

From the above relations, a solution of (\ref{coupleEig}) is
explicitly given by linear combinations (possibly $\BC_{2n}-$valued)
using monomials of the type
$$(x)_{+}^{(\alpha^+)}(x)_{-}^{(\alpha^-)},$$ with $|\alpha^+|=p$
and $|\alpha^-|=q$.

The crucial difference between the homogeneous polynomials
formulated in~\cite{BSS07} and the solutions of the coupled problem
(\ref{coupleEig}) is that they do not satisfy $$R_{p,q}(\eta
z,\tilde{\eta}z^\dag)=\eta^p(\tilde{\eta})^qR_{p,q}( z,z^\dag).$$
However, they form a Sheffer sequence of polynomials with respect to
the operators $\partial_h^{\pm j}$ and $x_jT_h^{\mp j}$, since
$(x)_{\mp}^{(\alpha)}$ are basic monomials (i.e.~$\partial_h^{\pm j}
x_k=\delta_{jk},~(\underline{x})_{\pm}^{({\bf 0})}={\bf 1}$ and $(\underline{{\bf 0}})_{\pm}^{(\alpha)}=0$.) and satisfy the
monomial principle
\begin{eqnarray*}
x_jT_h^{\mp j}(x)_{\mp}^{(\alpha)}=(x)_{\mp}^{(\alpha+ \vv_j)}, & \partial_h^{\pm j }(x)_{\mp}^{(\alpha)}=\alpha_j(x)_{\mp}^{(\alpha- \vv_j)} \\
\end{eqnarray*}
A further consequence of the above relations is the Rodrigues
formula
\begin{eqnarray*}
(x)_{\pm }^{(\alpha)}=(x_1T_{h}^{\pm 1})^{\alpha_1}(x_2T_{h}^{\pm 2})^{\alpha_2} \ldots (x_nT_{h}^{\pm n})^{\alpha_n}{\bf 1}.
\end{eqnarray*}
Moreover, $\left(z^p~{\bf 1}\right)((z^\dag)^q~{\bf 1})$ is a discrete homogeneous polynomial of degree $(p,q)$.

The above formulae together with the Weyl-Heisenberg character of
$\partial_h^{\pm j}$ and $x_j T_h^{\mp j}$ clearly suggests a
correspondence between the $\BR-$polynomial algebra generated by
$(x)_{\pm}^{(\alpha)}$ and the Bose algebra. In fact an isomorphism
between both algebras already tacitly exists \cite{BLR98}.

This is the starting point for constructing discrete versions of
hermitian Fischer decompositions in terms of discrete hermitian monogenic
homogeneous polynomials, i.e. polynomial solutions satisfying
(\ref{coupleEig}) and (\ref{hmon}) and, moreover, generate Hermite
polynomials as an Appel sequence associated with the orthogonally
shift-invariant Weierstrass operator~\cite{BLR98}. This will be one
of the main topics to be studied on the forthcoming
paper~\cite{nfaust08} from the umbral calculus point of view.

\subsection*{Acknowledgments}

We are grateful to Prof. E.~Forgy and I.~Kanamori for calling our attention to their papers~\cite{FS04} and \cite{KK04}.

During the preliminary redaction of this paper, the first author received the sad news that her cousin and best friend V\^ania Alexandre died on November 10, 2007. Therefore, this paper is also dedicated to her memory.


\end{document}